\documentclass{article}

\usepackage{t1enc}
\usepackage[latin2]{inputenc}

\usepackage{amssymb,amsmath}
\usepackage[slovene, croatian, slovak, english]{babel}

\newenvironment{proof}{\begin{trivlist}\item[]{\it
Proof.}}{\hfill$\square$\end{trivlist}}

\newtheorem{theorem}{Theorem}[section]

\newtheorem{conj}[theorem]{Conjecture}

\newcommand{\R}{\mathbb R}

\newcommand{\la}{\lambda}

\newcommand{\rk}{\mathrm {rk}\;}
\newcommand{\ifff}{if and only if }

\newcommand{\pd}{positive definite}
\newcommand{\haf}{\mathrm {haf}\;}
\newcommand{\pf}{\mathrm {pf}\;}
\newcommand{\per}{\mathrm {per}\;}

\begin{document}
\overfullrule=5pt

\date{}

\title{Pfaffians, hafnians and products of real linear functionals}

\author{ P\'eter E.\ Frenkel\footnote{ Partially supported by OTKA grants T 046365, K 61116 and NK 72523.}\\Alfr\'ed R\'enyi Institute of Mathematics\\ Hungarian Academy of Sciences\\ 
 P.O.B.\ 127,  1364 Budapest, Hungary\\{\tt frenkelp@renyi.hu}} 

\maketitle
\begin{abstract} 
We prove pfaffian and hafnian versions of Lieb's inequalities on determinants
and 
permanents of
positive
semi-definite matrices. We use the
hafnian inequality to improve the lower bound of R\'ev\'esz and Sarantopoulos
on the norm of a product of  linear functionals on a real Euclidean space (this subject is
sometimes called the `real linear polarization constant' problem).


\end{abstract} 


\section*{-1.  Introduction}
The contents of this  paper are as follows. In Section 0, we sketch one 
part of
the
historic background: classical inequalities on determinants
and 
permanents of
positive
semi-definite matrices. In Section~\ref{new}, we prove pfaffian and hafnian
versions of these inequalities, and we formulate Conjecture~\ref{conj}, another
hafnian inequality.
In Section~\ref{prod}, we apply the hafnian inequality of
Theorem~\ref{spec} to our main goal: improving 
the lower bound of R\'ev\'esz and Sarantopoulos
on the norm of a product of  linear functionals on a real Euclidean space (this subject is
sometimes called the `real linear polarization constant' problem, its history
is sketched at the end of the paper). This is achieved in
Theorem~\ref{polar}. We point out that Conjecture~\ref{conj} would be
sufficient to completely settle  the real linear polarization constant problem.

\section*{0.  Old inequalities on determinants 
and permanents}\label{old}
 Recall that the determinant and the permanent of an $n\times
n$ matrix $A=(a_{i,j})$ are defined by $$\det
A=\sum_{\pi\in\mathfrak S_n}(-1)^\pi\prod_{i=1}^na_{i,\pi(i)},\qquad\qquad\per A=\sum_{\pi\in\mathfrak S_n}\prod_{i=1}^na_{i,\pi(i)},$$
where $\mathfrak S_n$ is the symmetric group on $n$ elements. Throughout this
section, we assume that  $A$
is a positive semi-definite Hermitian $n\times n$ matrix (we
write $A\ge 0$). For such $A$,  Hadamard  proved that $$\det
A\le\prod_{i=1}^na_{i,i},$$ with equality if and only if  $A$ has a zero row or is a
diagonal matrix.
Fischer  generalized this to $$\det A\le\det A'\cdot\det A''$$ for 
\begin{equation}\label{bl}A=\left(\begin{matrix}
      A'&B\\B^*&A''\end{matrix}\right)\ge 0,\end{equation}
 with equality if and only if $\det
  A'\cdot\det A''\cdot B=0$. 

Concerning the permanent of a positive semi-definite matrix, Marcus [Mar1, Mar2] proved that
 \begin{equation}\label{M}\per A\ge \prod_{i=1}^na_{i,i},\end{equation}
 with equality if and only if  $A$ has a zero row or is a
diagonal matrix. Lieb [L] generalized this to \begin{equation}\label{L}
\per A\ge\per A'\cdot \per
A''\end{equation}
for  $A$ as in \eqref{bl},  with equality if and only if  $A$ has a zero row or $B=0$.
Moreover, he proved that in the polynomial $P(\lambda)$ of degree $n'$ (=size
of $A'$) defined by  
$$P(\lambda)=\per\left(\begin{matrix}
      \lambda A'&B\\B^*&A''\end{matrix}\right)=\sum_{t=0}^{n'}c_t\lambda ^t,$$
  all coefficients $c_t$ are real and non-negative. This is indeed a stronger
theorem since it implies $$\per A=P(1)=\sum_{t=0}^{n'} c_t\ge c_{n'}=\per A'\cdot\per A''.$$
\DJ okovi\'c [D, Mi] gave a simple proof of Lieb's inequalities, and showed also that
if $A'$ and $A''$ are positive definite then $c_{n'-t}=0$ if and only if all
subpermanents of $B$ of order $t$ vanish. Lieb [L] 
also states an analogous (and analogously provable) theorem for determinants:
for $A$ as in \eqref{bl}, let  $$D(\lambda)=\det\left(\begin{matrix}
      \lambda A'&B\\B^*&A''\end{matrix}\right)=\sum_{t=0}^{n'}d_t\lambda ^t.$$
If $\det A'\cdot\det A''=0$, then $D(\la)=0$. If $A'$ and $A''$ are positive
definite, then $(-1)^td_{n'-t}$  is
  positive for $t\le\rk B$ and is zero for $t>\rk B$.

\bigskip

\bf Remark. \rm In all of Lieb's inequalities mentioned above, 
the condition that the matrix $A$ is positive semi-definite can be replaced by
the weaker condition that the diagonal blocks $A'$ and $A''$ are positive semi-definite.  
The proof goes through virtually unchanged. Alternatively, this stronger form
of the inequalities can be easily deduced from the seemingly weaker form above.

\section{New inequalities on pfaffians and hafnians}\label{new}
For an $n\times n$ matrix $A=(a_{i,j})$ and subsets $S$,  $T$ of $N:=\{1,\dots, n\}$, we write $A_{S,T}:=(a_{i,j})_{i\in
  S, j\in T}.$ 
If $|T|=2t$ is even, we write
$$(-1)^T:=(-1)^{t+\sum_{j\in T}j}.$$ 

\subsection{Pfaffians} As far as the applications in Section~\ref{prod} are
concerned, this subsection may be skipped.

 Recall that the pfaffian  of a $2n\times
2n$ antisymmetric matrix $C=(c_{i,j})$ is  defined by $$\pf
C=\frac 1{n!2^n}\sum_{\pi\in\mathfrak S_{2n}}(-1)^\pi
c_{\pi(1),\pi(2)}\cdots c_{\pi(2n-1),\pi(2n)}.$$ We have
$\left(\pf C\right)^2=\det C$.

For antisymmetric $A$ and symmetric $B$, both of size $n\times n$,
we consider the polynomial $$(-1)^{\lfloor n/2\rfloor}
\pf \left(\begin{matrix} -\la A & B\\ -B &A
  \end{matrix}\right)=\sum_{t=0}^{\lfloor n/2\rfloor} p_t\la ^t.$$

\begin{theorem}\label{pf} Let $A$ and $B$ be  real
  $n\times n$ matrices with $A$ antisymmetric and $B$ symmetric. 
If $B$ is positive semi-definite, then $p_t\ge 0$
for all $t$. If $B$ is positive definite, then $p_t>0 $ for $t\le (\rk A)/2$
and  $p_t=0 $ for $t> (\rk A)/2$.
\end{theorem}

\begin{proof} If  $B=(b_{i,j})$ is positive semi-definite, then
there exist vectors 
$x_1$, \dots, $x_n$ in a real 
Euclidean space $V$ such that $(x_i,x_j)=b_{i,j}.$ Recall that in the exterior
tensor algebra $\bigwedge V$ a positive definite inner product (and the
corresponding Euclidean norm) is defined by
$$\left(\bigwedge v_i,\;\bigwedge w_j\right):=\det ((v_i,w_j)).$$
We have \begin{align*}
p_t=\sum_{|S|=2t}\sum_{|T|=2t}(-1)^S(-1)^T\pf A_{S,S}\cdot
\pf A_{T,T}\cdot\det B_{N\setminus S,N\setminus T }= \\
=\sum_{|S|=2t}\sum_{|T|=2t}\left((-1)^S\pf A_{S,S}\cdot \bigwedge_{i\not\in
    S}x_i
,\;  (-1)^T\pf A_{T,T}\cdot \bigwedge_{j\not\in
    T}x_j\right)=\\
=\left|\sum_{|S|=2t}(-1)^S\pf
  A_{S,S}\cdot\bigwedge_{i\not\in S}x_i\right|^2\ge 0.\end{align*}
Assume that $B$ is \pd. Then the vectors $x_i$ are linearly independent. It
follows that the tensors $\bigwedge_{i\not\in S}x_i$ are also linearly
independent as $S$ runs over the subsets of $N$. Thus $p_t=0$ if and only if  $\pf A_{S,S}=0$
for all $|S|=2t$, i.e., if and only if $2t>\rk A$.
\end{proof}

\begin{theorem}\label{pfspec} Let $A$ and $B$ be  real
  $n\times n$ matrices with $A$ antisymmetric and $B$ symmetric. 
Let $\la\ge 0$.
If $B$ is positive semi-definite, then   $$(-1)^{\lfloor n/2\rfloor}
\pf \left(\begin{matrix} -\la A & B\\ -B &A
  \end{matrix}\right)\ge\det B.$$
If $B$ is positive definite, then equality occurs \ifff    $\la A=0$.
\end{theorem}

\begin{proof} 
The left hand side is $$p_0+p_1\la+\dots+p_{\lfloor n/2\rfloor}
\la^{\lfloor n/2\rfloor}.$$ The
  right hand side is $p_0$.
\end{proof}

I am grateful to the anonymous referee of this paper for the idea of the
following alternative proof of Theorems~\ref{pf} and \ref{pfspec}.  We may
assume $B>0$, since every positive semi-definite matrix is a limit of positive
definite ones.  The matrix $B^{-1/2}AB^{-1/2}$ being real and antisymmetric,
there exists a unitary matrix $U$ such that $D:=U^{-1}B^{-1/2}AB^{-1/2}U$ is
diagonal with purely  imaginary eigenvalues $a_1\sqrt{-1}$, \dots,
$a_n\sqrt{-1}$.
The real multiset $\{a_1,\dots, a_n\}$ is invariant under $a\leftrightarrow
-a$.
We have \begin{align*}
\left(\sum p_t\la ^t\right)^2
=\det
\left(\begin{matrix} -\la A & B\\ -B &A
  \end{matrix}\right)
=\det
\left(\begin{matrix} -\la \sqrt BUDU^{-1} \sqrt B & B\\ -B &\sqrt BUDU^{-1}
    \sqrt B
  \end{matrix}\right)=\\
=\det\left(\left(\begin{matrix} \sqrt BU & 0\\ 0
      &\sqrt BU
  \end{matrix}\right)\left(\begin{matrix} -\la D & \bf 1 \\ -\bf 1  &D
  \end{matrix}\right)\left(\begin{matrix}U^{-1} \sqrt B & 0\\ 0 &U^{-1}
    \sqrt B
  \end{matrix}\right)\right)=\\
=\det\sqrt B^4\cdot\prod_{i=1}^n\det\left(\begin{matrix} -\la a_i\sqrt{-1} & 1\\ -1 &a_i\sqrt {-1}
  \end{matrix}\right)=\det B^2\cdot\prod_{i=1}^n(1+a_i^2\la).\end{align*}
Extracting square roots, and choosing the sign in accordance with $p_0=+\det
B$, we get $$\sum p_t\la ^t=(-1)^{\lfloor n/2\rfloor}
\pf \left(\begin{matrix} -\la A & B\\ -B &A
  \end{matrix}\right)=\det B\cdot\prod_{a_i>0}(1+a_i^2\la),$$ whence
both theorems immediately follow, since $\det B > 0$.

\subsection{Hafnians}
 Recall that the hafnian  of a $2n\times
2n$ symmetric matrix $C=(c_{i,j})$ is  defined by $$\haf
C=\frac 1{n!2^n}\sum_{\pi\in\mathfrak S_{2n}} c_{\pi(1),\pi(2)}\cdots 
c_{\pi(2n-1),\pi(2n)}.$$
For symmetric $A$ and $B$, both of size $n\times n$,
we consider the polynomial $$\haf \left(\begin{matrix} \la A & B\\ B &A
  \end{matrix}\right)=\sum_{t=0}^{\lfloor n/2\rfloor} 
h_t\la ^t.$$

\begin{theorem}\label{haf} Let $A$ and $B$  be symmetric  real
  $n\times n$ matrices. 
If $B$ is positive semi-definite, then $h_t\ge 0$
for all $t$. If $B$ is positive definite, then $h_t=0$ if and only if all
$2t\times 2t$ subhafnians of  $A$ vanish.
\end{theorem}

\begin{proof} If  $B=(b_{i,j})$ is positive semi-definite,  then
there exist vectors 
$x_1$, \dots, $x_n$ in a real 
Euclidean space $V$ such that $(x_i,x_j)=b_{i,j}.$ Recall [Mar1, Mar2, MN, Mi]
 that in the symmetric
tensor algebra $S V$ a positive definite inner product (and the
corresponding Euclidean norm) is defined by
$$\left(\prod v_i,\prod w_j\right):=\per ((v_i,w_j)).$$
We have \begin{align*}h_t=\sum_{|S|=2t}\sum_{|T|=2t}\haf A_{S,S}\cdot
\haf A_{T,T}\cdot\per B_{N\setminus S,N\setminus T }= \\
=\left|\sum_{|S|=2t}\haf
  A_{S,S}\cdot\prod_{i\not\in S}x_i\right|^2\ge 0.\end{align*}
Assume that $B$ is \pd. Then the vectors $x_i$ are linearly independent. It
follows that the tensors $\prod_{i\not\in S}x_i$ are also linearly
independent as $S$ runs over the subsets of $N$. Thus $h_t=0$ if and only if  $\haf A_{S,S}=0$
for all $|S|=2t$.
\end{proof}

\begin{theorem}\label{spec}
 Let $A$ and $B$  be symmetric  real
  $n\times n$ matrices. Let $\la\ge 0$. 
If $B$ is positive semi-definite, then
$$\haf\left(\begin{matrix} \la 
A&B\\ B&A\end{matrix}\right)\ge\per B.$$ If $B$ is
  positive definite, then equality occurs if and only if $A$ is a diagonal
  matrix
or $\la=0$.
\end{theorem}

\begin{proof} The left hand side is $$h_0+h_1\la+\dots+h_{\lfloor n/2\rfloor}
\la ^{\lfloor n/2\rfloor}.$$ The
  right hand side is $h_0$.
\end{proof}

Setting $A=B$ and $\la =1$, and combining with Marcus's inequality \eqref{M},
we arrive at case $p=1$ of 

\begin{conj}\label{conj} If $A=(a_{i,j})$ is a positive semi-definite symmetric real $n\times n$
  matrix, then the hafnian of the $2pn\times 2pn$ 
matrix consisting of $2p\times
  2p$ blocks $A$ is at least ${(2p-1)!!}^n\prod a_{i,i}^p,$ with equality \ifff
  $A$ has  a zero row or is a  diagonal matrix.
\end{conj}

\section{Products of real linear functionals}\label{prod}
In this section,  we apply Theorem~\ref{spec} to products of jointly normal
random variables and then to products of real linear
functionals, which was the main motivation for this work.
The ideas in this section are analogous to those that
 Arias-de-Reyna [A] used in
the complex case.

Let $\xi_1$, \dots, $\xi_d$
 denote  independent random  variables with standard
Gaussian distribution, i.e., with
joint density function
$(2\pi)^{-d/2}\exp({-|\xi|^2/2})$, where $|\xi|^2=\sum\xi_k^2.$
We write $Ef(\xi)$ for the expectation 
of a function $f=f(\xi)=f(\xi_1,\dots, \xi_d)$. 
Recall that 
$$
E\xi_k^{2p}=(2p-1)!!=(2p-1)(2p-3)\cdots 3\cdot 1$$ for $k=1,\dots, d$ (easy inductive proof via integration
by parts), and thus $$E\prod_{k=1}^d 
\xi_k^{2p_k}=\prod_{k=1}^d (2p_k-1)!!.$$ 

On $\R^d$, we write 
$(\cdot,\cdot)$ 
for the standard Euclidean inner product.
We recall the well-known [B2, G, S, Z]

\bigskip\noindent {\bf Wick formula.}\; \it 
Let $x_1$, \dots, $x_n$ be vectors in $\R^d$ with Gram matrix $A=((x_i,x_j)).$ Then \begin{equation}\label{Sigma}
E\prod_{i=1}^n (x_i,\xi)=\haf A. 
\end{equation}
\rm
(For odd $n$, we define $\haf A=0$.)
\begin{proof}
Both sides are multilinear in the $x_i$, so we may assume 
that each $x_i$ is an element of the standard orthonormal basis $e_1$, \dots,
$e_d$. If there is an $e_k$ 
that occurs an odd number of   times among the $x_i$, then 
both sides are zero. If each  $e_k$ 
occurs $2p_k$ times, then the left hand side is $E\prod_{k=1}^d 
\xi_k^{2p_k}$, and the right hand side is
$\prod_{k=1}^d (2p_k-1)!!$, which are equal.
\end{proof}

The following  theorems are easy corollaries of  Theorem~\ref{spec}
together with the Wick formula~\eqref{Sigma} and Marcus's theorem \eqref{M}. 

\begin{theorem}\label{mom} If $X_1$, \dots, $X_n$ are jointly normal random
  variables with zero expectation, then
$$E\left( X_1^2\cdots X_n^2\right)\ge EX_1^2\cdots EX_n^2.$$
Equality holds  if and only if they are independent or at least one of them is almost
surely zero.
\end{theorem}

\begin{proof} The variables can be written as $X_i=(x_i,\xi)$ with $\xi$  of standard
normal distribution and the $x_i$ constant vectors with a positive
semi-definite Gram matrix $A=(a_{i,j})=((x_i,x_j))$. 
Then \begin{align*}E\prod _{i=1}^n X_i^2=E\prod _{i=1}^n (x_i,\xi)^2=\\
=\haf\left(\begin{matrix} A&A\\
      A&A\end{matrix}\right)\ge\per
  A\ge\prod_{i=1}^na_{i,i}=\\
=\prod_{i=1}^nE(x_i,\xi)^2=\prod_{i=1}^nEX_i^2,\end{align*} with equality if and
  only if $A$
  is a diagonal matrix or has a zero row, i.e., the $x_i$ are pairwise
  orthogonal or at least one of them is zero.
\end{proof}

The generalization of Theorem~\ref{mom} to an arbitrary even exponent $2p$ is
equivalent to Conjecture~\ref{conj}.

\newcommand{\ii}{\'{\i}} 

\begin{theorem}\label{average} For any $x_1,\dots, x_n\in\R^d$, $|x_i|=1$, the average of $\prod
(x_i,\xi)^2$ on the unit sphere $\{\xi\in\R^d\;:\; |\xi|=1\}$ is at least
$$\frac{\Gamma
  (d/2)}{2^n\Gamma (d/2+n)}=\frac {(d-2)!!}{(d+2n-2)!!}=\frac 1{d(d+2)(d+4)\dots
(d+2n-2)},$$ with equality  if and only if the vectors $x_i$
are pairwise orthogonal. 
\end{theorem}

\begin{proof} The average on the unit sphere is the constant in the theorem
  times the expectation w.r.t.\ the standard Gaussian measure (see e.g.\ [B1]). By
  Theorem~\ref{mom}, the latter expectation  is minimal \ifff the $x_i$ are pairwise orthogonal,
  in which case it is 1.
\end{proof}

\begin{theorem}\label{polar} For real linear functionals $f_i$ on a real Euclidean space,
  $$||f_1\cdots f_n||\ge \frac{||f_1||\cdots
  ||f_n||}{\sqrt{n(n+2)(n+4)\cdots (3n-2)}}.$$
\end{theorem}
Here $||\cdot ||$ means supremum of the absolute value on the unit sphere.
In the infinite-dimensional case, functionals with infinite norm may be
allowed. Then the convention $0\cdot\infty=0$ should be used on the right hand
side.

\begin{proof} We may assume that the space is $\R^d$ with $d\le n$,  and the functionals are given by
  $f_i(\xi)=(x_i,\xi)$ with $||f_i||=|x_i|=1.$ Then  $||f_1\cdots f_n||^2$
is at least the average of $\prod f_i^2(\xi)=\prod (x_i,\xi)^2$ on the unit sphere, which by
Theorem~\ref{average} and $d\le n$
is at least $1/(n(n+2)(n+4)\cdots (3n-2)).$
\end{proof} 

It is   an unsolved problem, raised by Ben\ii tez, Sarantopoulos and Tonge
[BST] (1998),
whether Theorem~\ref{polar} is  true
with $n^n$ under the square root sign
in the denominator on the right hand side. This is called the `real linear
polarization constant' problem. In the complex case, the affirmative answer
was proved by Arias-de-Reyna [A] in 1998, based on  the complex analog of the
Wick formula [A, B2, G] and on Lieb's inequality \eqref{L}.\footnote{The
 referee of the present paper called my attention to the fact that
  Arias-de-Reyna used only the special case of \eqref{L} where the matrix $A'$
  is
of rank 1. This is much simpler than \eqref{L} in general, it can be proved
essentially
by the argument Marcus used in [Mar1, Mar2] to prove the 
 even more special case $n'=1$, which still implies
inequality~\eqref{M}.} Keith Ball [Ball] gave another proof of the affirmative
answer in the complex case by solving the
complex plank problem.

In the real case, the affirmative answer 
for $n\le 5$ was proved by Pappas and R\'ev\'esz [PR] in 2004.
For general $n$, the best estimate  known before the present 
paper was that of R\'ev\'esz and Sarantopoulos
[RS] (2004), based on results of [MST],
with $(2n)^{n}/4$ under the square root sign. 
See [Mat1, Mat2, MM, R] for accounts on this and
related questions.
Note that \begin{align*}
{n(n+2)(n+4)\cdots (3n-2)}=\\
=\exp\left(\log n+\log (n+2)+\log (n+4)+\dots+\log
(3n-2)\right)<\\
<\exp\left({\frac12\int_n^{3n}\log u\cdot {\rm d}u}\right)=\\
=\exp\left(
 [u(\log u
-1)]_n^{3n}/2\right)=\exp ((3n\log {3n}-3n
-n\log n+n)/2)=\\
=\exp {\frac {n(2\log n+3\log 3 -2)}2}
=\left(\frac{3\sqrt 3}en\right)^n,\end{align*}
 and $3\sqrt
3/e<3\cdot 1.8/2.7=2$, so Theorem~\ref{polar} is an improvement.
Note also that the statement with $n^n$ 
under  the square root sign
would follow from
Conjecture~\ref{conj}.

\section*{Acknowledgements}
I am grateful to P\'eter Major, M\'at\'e Matolcsi and  Szil\'ard R\'ev\'esz 
for
helpful discussions, and to  the  anonymous referee for useful comments.

\section*{References}
\noindent
[A] J.\ Arias-de-Reyna, Gaussian variables, polynomials and permanents, Lin.\
Alg.\ Appl.\ 285 (1998), 107--114.

\bigskip\noindent
[Ball] K.\ M.\ Ball, The complex plank problem, Bull.\ London.\ Math.\ Soc.\
33 (2001), 433--442.

\bigskip\noindent
[B1] A.\ Barvinok, Estimating $L^\infty$ norms by $L^{2k}$ norms for functions
on orbits, Found.\ Comput.\ 
Math.\ 2 (2002), 393--412.

\bigskip\noindent
[B2] A.\ Barvinok, Integration and optimization of multivariate polynomials by
restriction onto a random subspace, arXiv preprint:
math.OC/0502298

\bigskip\noindent
[BST] C.\ Ben\ii tez, Y.\ Sarantopoulos, A.\ Tonge, Lower bounds for norms of
products of polynomials, Math.\ Proc.\ Camb.\ Phil.\ Soc.\ 124 (1998),
395--408.

\bigskip\noindent
[D] D.\ \v Z.\ \DJ okovi\'c, Simple proof of a
theorem on permanents, Glasgow
Math.\ J.\ 10 (1969), 52--54.

\bigskip\noindent
[G] L.\ Gurvits, Classical complexity and quantum entanglement,  J.\ Comput.\
 System Sci.\  69  (2004),  no. 3, 448--484.

\bigskip\noindent
[L] E.\ H.\ Lieb, Proofs of some conjectures on permanents, J.\ Math.\ Mech.\
16 (1966), 127--134.

\bigskip\noindent
[Mar1]  M.\ Marcus, The permanent analogue of the Hadamard determinant theorem,
Bull.\ Amer.\ Math.\ Soc.\ 69 (1963), 494--496.

\bigskip\noindent
[Mar2] M.\ Marcus, The Hadamard theorem for permanents, Proc.\ Amer.\ Math.\
Soc.\ 15 (1964), 967--973.

\bigskip\noindent
[MN] M.\ Marcus, M.\ Newman, The permanent function as an inner product,
Bull.\ Amer.\ Math.\ Soc.\ 67 (1961), 223--224.

\bigskip\noindent
[Mat1] M.\ Matolcsi, A geometric estimate on the norm of  product of functionals,
Lin.\ Alg.\ Appl.\ 405 (2005), 304--310. 

\bigskip\noindent
[Mat2] M.\ Matolcsi, The linear polarization constant of $\R^n$, Acta Math.\
Hungar.\ 108 (2005), no.\ 1-2, 129--136.

\bigskip\noindent
[MM] M.\ Matolcsi, G.\ A.\ Mu\~noz, On the real linear polarization constant
problem, Math.\ Inequal.\ Appl.\ 9 (2006), no.\ 3, 485--494.

\bigskip\noindent
[Mi] H.\ Minc, Permanents, Encyclopedia of Mathematics and its Applications,
Add\-is\-on-Wesley, 1978

\bigskip\noindent
[MST]
 G.\ A.\ Mu\~noz, Y.\ Sarantopoulos, A.\ Tonge, Complexifications of real
 Banach spaces, polynomials and multilinear maps, Studia Math.\ 134 (1999),
 no.\ 1, 1--33.

\bigskip\noindent
[PR] A.\ Pappas, Sz.\ R\'ev\'esz, Linear polarization constants..., J.\ Math.\
Anal.\ Appl. 300 (2004), 129--146.

\bigskip\noindent
[R] Sz.\ Gy.\
R\'ev\'esz, Inequalities for multivariate polynomials, Annals of the Marie
Curie Fellowships 4 (2006), {\tt http:/\!/www.mariecurie.org/annals/},
arXiv preprint: math.CA/0703387

\bigskip\noindent
[RS] Sz.\ Gy.\
R\'ev\'esz, Y.\ Sarantopoulos, Plank problems, polarization and Chebyshev
constants, J.\ Korean Math.\ Soc.\ 41 (2004) 157--174.

\bigskip\noindent
[S] B.\ Simon,
The P$(\phi)_2$ Euclidean (Quantum) Field Theory, Princeton Series in Physics,
Princeton University Press, 1974

\bigskip\noindent
[Z] A.\ Zvonkin, Matrix integrals and map enumeration: an accesible
introduction, Combinatorics and physics (Marseille, 1995), Math.\ Comput.\
Modelling 26 (1997), 281--304.
\end{document}